\documentclass[11pt]{article}
\usepackage{amsfonts,amsmath,amssymb}
\usepackage{algorithm2e,framed}

\newcommand{\Expect}[1]{\mbox{}{\bf{E}}\left[#1\right]}

\newcommand{\FNorm }[1]{\mbox{}\left\|#1\right\|_F  }
\newcommand{\FNormS}[1]{\mbox{}\left\|#1\right\|_F^2}

\newcommand{\TNorm }[1]{\mbox{}\left\|#1\right\|_2  }
\newcommand{\TNormS}[1]{\mbox{}\left\|#1\right\|_2^2}

\newcommand{\XNorm }[1]{\mbox{}\left\|#1\right\|_{\xi}  }

\newcommand{\VTTNormS}[1]{\mbox{}\left\|#1\right\|_2^2}

\newcommand{\abs }[1]{\left|#1\right|}

\long\def\killtext#1{}
\newtheorem{problem}{Problem}
\newtheorem{definition}{Definition}

\newtheorem{lemma}{Lemma}
\newtheorem{theorem}{Theorem}

\newenvironment{Proof}{\noindent {\em Proof:}}{\\\hspace*{\fill}\mbox{$\diamond$}}

\newcommand{\nnz}[1]{{\bf nnz}{\left(#1\right)}}

\newcommand{\argmin}{\text{argmin}}

\begin{document}

\title{
Effective Resistances, Statistical Leverage, and \\
Applications to Linear Equation Solving
}

\author{
Petros Drineas
\thanks{
Department of Computer Science,
Rensselaer Polytechnic Institute,
Troy, NY, drinep@cs.rpi.edu.
}
\and
Michael W. Mahoney
\thanks{
Department of Mathematics,
Stanford University,
Stanford, CA,
mmahoney@cs.stanford.edu.
}
}

\date{}
\maketitle

\begin{abstract} 
Recent work in theoretical computer science and scientific computing has
focused on nearly-linear-time algorithms for solving systems of linear
equations.
While introducing several novel theoretical perspectives, this work has
yet to lead to practical algorithms.
In an effort to bridge this gap, we describe in this paper two related
results.
Our first and main result is a simple algorithm to approximate the solution
to a set of linear equations defined by a Laplacian (for a graph $G$ with
$n$ nodes and $m \le n^2$ edges) constraint matrix.
The algorithm is a non-recursive algorithm; even though it runs in
$O(n^2 \cdot \mbox{polylog}(n))$ time rather than $O(m \cdot \mbox{polylog}(n))$ time 
(given an oracle for
the so-called statistical leverage scores), it is extremely simple; and it
can be used to compute an approximate solution with a direct solver.
In light of this result, our second result is a straightforward connection
between the concept of graph resistance (which has proven useful in recent
algorithms for linear equation solvers) and the concept of statistical
leverage (which has proven useful in numerically-implementable randomized
algorithms for large matrix problems and which has a natural data-analytic
interpretation).
\end{abstract}


\section{Introduction}
\label{sec:intro}

The problem of approximating the solution to a set of linear equations
defined by a Laplacian constraint matrix has been of interest recently due
to a series of remarkable papers by Spielman and
Teng~\cite{ST08b_TR,ST08a_TR,ST06a_TR}.
(This work builds on ideas originally introduced by Vaidya and developed
by others~\cite{BH03,BCHT04,BGHNT06}.%
\footnote{Briefly, recall that to solve a system of linear equations,
$Ax=b$, one can use either direct methods or iterative
methods~\cite{GVL96,TrefethenBau97}.
Iterative methods, such as Chebyshev or Conjugate Gradients, compute
successively better approximations to $x$ by performing successive
matrix-vector multiplications.
The number of iterations typically depends on the condition number
$\kappa(A)$ of $A$, where $\kappa(A)=\lambda_{max}(A)/\lambda_{min}(A)$ is
the ratio of the extreme (nontrivial) eigenvalues of $A$, via a
multiplicative factor of $\sqrt{\kappa(A)}$.
Preconditioning refers to a class of methods to solve $B^{-1}Ax=B^{-1}b$,
where the preconditioning matrix $B$ is chosen such that $\kappa(B^{-1}A)$
is small and such that it is easy to solve for $Bz=c$.
Vaidya introduced the idea of using combinatorial methods to precondition
Laplacians of graphs with Laplacians of their subgraphs.
It is known that if one wants to precondition any symmetric diagonally
dominant matrix, then it suffices to find a preconditioner for a related
Laplacian matrix~\cite{BGHNT06};
and, moreover, that preconditioning matrices that arise in many applications
can be reduced to the problem of preconditioning diagonally dominant
matrices~\cite{BHV08}.
Vaidya's methods have been extended~\cite{BGHNT06,BH03,BCHT04}, and they
were used by Spielman and Teng to approximate the solution to diagonally
dominant linear systems in time that is ``nearly-linear'' in the number of
nonzero entries in their defining matrices~\cite{ST06a_TR}.})
While introducing several novel theoretical perspectives, this work on
``nearly-linear-time'' algorithms has yet to lead to practical algorithms.
In this paper, we describe two related results in an effort to bridge this
theory-practice gap.

Our first and main result, to be described in 
Section~\ref{sxn:LinearEqnSolving}, is a simple algorithm for computing an 
approximate solution to a set of linear equations defined by a Laplacian 
constraint matrix. The 
simplicity of the algorithm permits us to identify a simple connection (that 
to our knowledge has been overlooked) with other recent work in the theory 
and (numerical and data) application of randomized algorithms for matrix 
problems. 
Thus, our second result, to be described in Section~\ref{sxn:graph_resistances}, is to identify and discuss the connection between the concept of \emph{statistical leverage} and the concept of \emph{graph resistance}. The latter concept has a long history in spectral graph theory~\cite{Chung:1997}, and recently it has proven useful in algorithms for linear equation solvers~\cite{SS08a_STOC,BSS09a_STOC}. The former concept also has a long history, but in statistics and diagnostic data analysis~\cite{ChatterjeeHadi88}. Moreover, recently, it has been demonstrated to be the key structural quantity to understand in order to bridge the theory-practice gap between theoretical work on randomized algorithms for large matrices and applications (both numerical-implementation and data-analysis applications) of this ``randomized matrix algorithm''
paradigm~\cite{Paschou07b,CUR_PNAS,DMMS07_FastL2_TRv3,RT08,AMT09_DRAFT}.

\subsection{Laplacian matrices}\label{sxn:LaplacianDef}

Consider a graph $G=(V,E)$ with $n$ vertices and $m$ weighted, undirected 
edges. 
We will assume that all the weights are positive. 
Then, we can construct the so-called Laplacian matrix 
$L \in \mathbb{R}^{n \times n}$ of $G$. Let $w_{ij} \geq 0$ denote the 
weight of the edge joining vertices $i$ and $j$; clearly $w_{ij}=0$ if no 
such edge exists. In the most common definition of the Laplacian matrix $L$, 
the off-diagonal entries of $L$ ($L_{ij}$, $i \neq j$) are set to $-w_{ij}$, 
while the diagonal entries $L_{ii}$ (for all $i = 1,\ldots,n$) are equal to 
the ``weighted degree'' of vertex $i$, i.e., $L_{ii} = \sum_{j=1}^n w_{ij}$. 
By definition, $L$ is a symmetric matrix of rank at most $n-1$, since the 
all-ones vector is clearly in the null space of $L$.

A somewhat less common definition of the Laplacian matrix follows from the 
so-called edge-incidence matrix of the graph $G$. 
Let $B \in \mathbb{R}^{m \times n}$ denote the edge-incidence matrix of the 
undirected graph $G$, constructed as follows: each row of $B$ corresponds to 
an edge of $G$; assuming that an (arbitrarily-oriented) edge of $G$ starts 
at vertex $i$ and ends at vertex $j$, the $i$-th entry in the corresponding 
row of $B$ is set to $+1$, the $j$-th entry is set to $-1$, and the 
remaining entries are all set to $0$. 
Thus, $B$ has two non-zero entries per row for a total of $2m$ non-zero 
entries. 
Also, let $W \in \mathbb{R}^{m \times m}$ be a diagonal matrix containing 
the edge weights (in the same order as they appear in $B$). 
Then, it is well-known that
$$L = B^T W B.$$
The above definition makes it obvious that $L$ is a symmetric 
positive-semidefinite matrix. 

Note that given a Laplacian matrix $L \in \mathbb{R}^{n \times n}$ corresponding to an undirected, weighted graph $G$ with $m$ edges and positive edge weights, we can immediately derive $B$ and $W$. 
That is, by considering the $m$ non-zero entries $L_{ij}$ with $i<j$, since 
each such entry corresponds to an edge joining vertices $i$ and $j$ of 
weight $w_{ij}=-L_{ij}$, we can immediately construct $B$ and $W$.

\subsection{An overview of the problem}

Given a Laplacian matrix $L$ corresponding to an underlying graph $G=(V,E)$ with $n$ vertices and $m$ (positively) weighted, undirected edges, consider the following regression problem which was addressed by Spielman and Teng~\cite{ST08b_TR,ST08a_TR,ST06a_TR}.
\begin{problem}
\label{prob:laplacian}
\noindent\textsc{[Least-squares approximation with Laplacian constraints]}
Given as input a Laplacian matrix $L \in \mathbb{R}^{n \times n}$ as described above and a target vector $b \in \mathbb{R}^n$, compute
\begin{equation*}
\arg \min_{x \in \mathbb{R}^n} \TNorm{Lx-b}.
\end{equation*}
The minimal $\ell_2$-norm solution vector $x_{opt}$ to the above problem is equal to
\begin{equation}\label{eqn:LaplacianOptVec}
x_{opt} = L^{\dagger}b,
\end{equation}
where $L^{\dagger}$ corresponds to the Moore-Penrose generalized inverse. 
\end{problem}
This formulation is a generalization of the standard problem of solving a system of linear equations of the form $Lx=b$, in order to better handle the rank-deficiency of $L$. We chose this formulation since our algorithm will make no assumptions on the rank of $L$. In addition, this formulation will make the comparison with related work on randomized algorithms for matrix problems (see Section~\ref{sxn:graph_resistances:usefulness}) immediate.

In this setting, Spielman and Teng~\cite{ST06a_TR} provided a randomized, 
relative-error approximation algorithm for Problem~\ref{prob:laplacian}. 
The running time of their algorithm is $O\left(\nnz{A} \log^{c_1} n\right)$, 
where $\nnz{A}$ represents the number of non-zero elements of the matrix 
$A$, or equivalently the number of edges in the graph $G$, and $c_1$ is a 
small constant. 
The first step of this algorithm corresponds to performing ``spectral graph sparsification,'' thereby keeping a small number of edges from $G$, and thus creating a much sparser Laplacian matrix $\tilde{L}$. The second step of this algorithm involves using this sparse matrix $\tilde{L}$ as an efficient preconditioner to solve Problem~\ref{prob:laplacian} approximately. In order to achieve high
precision, this is done in a recursive manner.

While~\cite{ST06a_TR} is a major theoretical breakthrough, its applicability is currently hindered by its sheer complexity. In an effort to bridge the gap between theory and practice, recent work of
Spielman and Srivastava~\cite{SS08a_STOC} proposed a much simpler algorithm for the graph sparsification step of~\cite{ST06a_TR}, by arguing that randomly sampling edges from the graph $G$ with probabilities proportional to the so-called \textit{effective resistances} (see Section~\ref{sxn:definitions} for definitions) of the edges provides a sparse Laplacian matrix $\tilde{L}$ satisfying the
desired properties. On the negative side, in order to approximate the effective resistances of the edges of $G$ efficiently, the Spielman-Srivastava algorithm performs $O(\log n)$ calls to the Spielman-Teng solver, severely hindering its applicability~\cite{SS08a_STOC}. We should also note that Batson, Spielman, and Srivastava~\cite{BSS09a_STOC} provided a more expensive algorithm for finding even sparser spectral sparsifiers.

Note that the work of Spielman and Teng also addresses a much broader class 
of matrices, the so-called $SDDM_0$ class, which can be reduced to the 
Laplacian case. 
This reduction is described in detail in~\cite{BHV08}. 
For simplicity of presentation, here we will only focus on Laplacian matrices.

\subsection{Solving systems of linear equations with Laplacian matrices}

Our main result in this paper is a simple algorithm to compute an approximate solution to Problem~\ref{prob:laplacian}. As with previous algorithms, the first phase will sparsify the input graph, and the second phase will solve the problem on the sparsified graph. Our main algorithm will be described in detail in Section~\ref{sxn:LinearEqnSolving}. Briefly, in the first phase we will compute a nonuniform sampling probability distribution that depends on the so-called statistical leverage scores~\cite{DMM08_CURtheory_JRNL,CUR_PNAS} associated with the weighted edge-incidence matrix of the input graph. We will then sample a ``small'' number of edges according to that distribution to construct a
sparsified Laplacian matrix $\tilde{L}$, having $O\left(\frac{n}{\epsilon} \log\frac{n}{\epsilon}\right)$ non-zero entries. Then, in the second phase we will solve the sparsified problem
\begin{equation}\label{eqn:LaplacianLS_Sample}
\arg \min_{x \in \mathbb{R}^n} \TNorm{\tilde{L}x-b}
\end{equation}
to get the vector $\tilde{x}_{opt} = L^{\dagger}b$. The resulting vector $\tilde{x}_{opt}$ satisfies (with constant probability)
\begin{equation}
\label{eqn:relerror}
\left\|x_{opt} - \tilde{x}_{opt} \right\|_L
   \leq \epsilon \left\|x_{opt}\right\|_L.
\end{equation}
Recall that the ``energy norm'' $\|x\|_L$ for any vector $x\in \mathbb{R}^n$ 
and any matrix $L \in \mathbb{R}^{n \times n}$ is equal to $x^T L x$. 
Given the sparsified Laplacian $\tilde{L}$, this second phase will use the 
conjugate gradient method as a direct solver~\cite{TrefethenBau97} to solve 
the sparse least-squares problem of eqn.~(\ref{eqn:LaplacianLS_Sample}), and 
thus it will take $O\left(\frac{n^2}{\epsilon} \log\frac{n}{\epsilon}\right)$ 
time. 
For dense graphs, this matches the running time of the Spielman-Teng 
algorithm, while for sparse graphs the Spielman-Teng algorithm is still 
faster.

The question of computing the statistical leverage scores (either exactly
or approximately) is a subtle one, and it is related to the theory-practice
disconnect---both for this problem, as well as for other problems to which
randomized matrix algorithms have been applied.
Thus, we will discuss this topic in greater detail in
Section~\ref{sxn:graph_resistances:usefulness}.
Briefly, $O(mn^2)$ time certainly suffices to compute them with standard
methods; theoretically, they can be computed in $O(m \log^{c_1} n)$ time, 
for some small constant $c_1$; and they can be efficiently approximated in 
the presence of certain resource constraints.


\section{An algorithm for solving systems of linear equations}
\label{sxn:LinearEqnSolving}

In this section, we will describe our main algorithm to approximate the minimal $\ell_2$-norm solution vector $x_{opt}$ of the least-squares approximation problem with Laplacian constraint matrix (Problem~\ref{prob:laplacian}). Then, we will state and prove our main quality-of-approximation theorem and discuss the running time of the proposed algorithm.

\subsection{Our main algorithm}

\begin{algorithm}[h]
\begin{framed}

\textbf{Input:} Laplacian matrix $L \in \mathbb{R}^{n \times n}$, corresponding to a graph $G$ with $n$ vertices and $m$ (positively) weighted edges, $b \in \mathbb{R}^n$, and accuracy parameter $\epsilon \in (0,1)$.

\vspace{0.05in}

\textbf{Output:}
$\tilde{x}_{opt} \in \mathbb{R}^n$.

\begin{enumerate}

\item Compute the edge-incidence matrix $B \in \mathbb{R}^{m \times n}$ and the diagonal edge-weight matrix $W \in \mathbb{R}^{m \times m}$ (see Section~\ref{sxn:LaplacianDef}).

\item Let $\Phi = W^{1/2}B \in \mathbb{R}^{m \times n}$.

\item Compute a set of probabilities $p_i$ (for all $i=1\ldots m$) such that $\sum_{i=1}^m p_i =1$ and
\begin{equation}\label{eqn:probabilities}
p_i \geq \frac{\beta \TNormS{\left(U_{\Phi}\right)_{(i)}}}{\FNormS{U_{\Phi}}}
\end{equation}
for some $\beta \in (0,1]$. ($U_{\Phi}$ is an orthogonal basis for the column space of $\Phi$ and $\left(U_{\Phi}\right)_{(i)}$ is the $i$-th row of $U_{\Phi}$.)
\item Set $r= \frac{72c_0^2n}{\beta \epsilon}\log\left(\frac{36c_0^2n}{\beta \epsilon}\right),$ where $c_0$ is the unspecified constant of Theorem~\ref{thm:theorem7correct}.

\item Initialize $S \in \mathbb{R}^{m \times r}$ to be an all-zeros matrix.

\item \textbf{For } $t=1,\ldots,r$ \textbf{do}

\begin{itemize}

\item Pick $i_t \in 1\ldots m$, where $\mbox{Prob}\left(i_t = i\right) = p_i$\;
\item $S_{i_t t} = 1/\sqrt{rp_{i_t}}$\;

\end{itemize}

\item Compute $\tilde{L} = \left(B^T W^{1/2}S\right) \left(S^TW^{1/2}B\right) \in \mathbb{R}^{n \times n}$.

\item Return $\tilde{x}_{opt} = \tilde{L}^{\dagger}b$.
\end{enumerate}

\end{framed}
\caption{Approximating the minimal $\ell_2$-norm solution of least-squares problems with Laplacian constraint matrices.} \label{alg:alg_sample_fast}
\end{algorithm}

Algorithm~\ref{alg:alg_sample_fast} takes as input an $n \times n$ Laplacian matrix $L$ (corresponding to a graph $G$ with $n$ vertices and $m$ positively weighted, undirected edges) and constructs an $n \times n$ sparsified Laplacian matrix $\tilde{L}$. Finally, it computes the minimal $\ell_2$-norm solution vector $\tilde{x}_{opt}$ of the sparsified problem with a direct solver.

In more detail, the algorithm first computes the edge incidence matrix $B$ and the corresponding diagonal weight matrix $W$, as described in Section~\ref{sxn:LaplacianDef}. Then, it computes a set of probabilities $p_1,p_2,\ldots,p_m$ such that the $i$-th edge of the graph, i.e., the $i$-th row of the edge-incidence matrix $B$ and the corresponding weight $W_{ii}$, will be retained with probability proportional to $p_i$. These probabilities satisfy eqn.~(\ref{eqn:probabilities}) and depend on the so-called statistical leverage scores of the matrix $W^{1/2}B$. 
As we will observe in Section~\ref{sxn:connection}, these scores are 
proportional to the effective resistances of the edges of graph $G$. 
The parameter $\beta$ at Step 3 of Algorithm~\ref{alg:alg_sample_fast} 
facilitates the use of approximate (as opposed to exact) probabilities and 
will be further discussed below.
It is worth noting that computing the aforementioned probabilities $p_i$ 
exactly ($\beta = 1$) necessitates $O(mn^2)$ time, which is prohibitive for 
the proposed application.%
\footnote{Indeed, one goal of this work is to focus further research towards 
efficient---either provably accurate or heuristic---algorithms to 
approximate these leverage scores in various settings, thereby leading to 
faster algorithms for this and related problems.}

After setting the sparsity parameter $r$ to an appropriate value that guarantees a relative-error approximation to the optimal solution at Step 4, exactly $r$ edges of $G$ are sampled (Step 6) with respect to the computed probabilities. The weights of the retained edges are rescaled (Step 6) and the induced Laplacian $\tilde{L}$ corresponding to the sparsified graph is formed. Note that $\tilde{L} \in \mathbb{R}^{n \times n}$ has at most $n+2r$ non-zero entries, since its underlying sparsified graph has at most $r$ edges. Then, the sparsified problem
\begin{equation}
\arg \min_{x \in \mathbb{R}^n} \TNorm{\tilde{L}x-b}
\end{equation}
is solved in order to return the minimal $\ell_2$-norm solution $\tilde{x}_{opt}=\tilde{L}^{\dagger}b$. The computational savings emerge since the sparsified problem can be solved efficiently using, for example, conjugate-gradient-type methods as direct solvers. The running time of such methods with input $\tilde{L}$ and $b$ is $O\left(n\left(n+2r\right)\right)$, where $n+2r$ is the number of non-zero entries in $\tilde{L}$.

\subsection{Approximation accuracy}

The following theorem is our main quality-of-approximation result for Algorithm~\ref{alg:alg_sample_fast}.
\begin{theorem}\label{thm:main}
Given Laplacian matrix $L \in \mathbb{R}^{n \times n}$ (corresponding to a graph $G$ with $n$ vertices and $m$ positively weighted edges) and a vector $b \in \mathbb{R}^n$, let $x_{opt} \in \mathbb{R}^n$ be the solution vector of eqn.~(\ref{eqn:LaplacianOptVec}). If $\tilde{x}_{opt} \in \mathbb{R}^n$ is the output of Algorithm~\ref{alg:alg_sample_fast} for some choice of the accuracy parameter $\epsilon \in (0,1)$, then, with probability at least 2/3,
\begin{equation*}
\left\|x_{opt} - \tilde{x}_{opt} \right\|_L \leq \epsilon \left\|x_{opt}\right\|_L.
\end{equation*}
\end{theorem}

\begin{Proof}
By definition, $\left\|x_{opt} - \tilde{x}_{opt} \right\|_L  = \left(x_{opt} - \tilde{x}_{opt}\right)^T L \left(x_{opt} - \tilde{x}_{opt}\right)$. Recall that $L = B^TWB$, where $B \in \mathbb{R}^{m \times n}$ and $W \in \mathbb{R}^{m \times m}$ are the edge-incidence and the diagonal weight matrix respectively (see Section~\ref{sxn:LaplacianDef}). Also recall that the diagonal entries of $W$ are positive and thus $W^{1/2}$ is well-defined. Then,
\begin{eqnarray}
\nonumber \left\|x_{opt} - \tilde{x}_{opt} \right\|_L  &=& \left(x_{opt} - \tilde{x}_{opt}\right)^T B^TWB \left(x_{opt} - \tilde{x}_{opt}\right)\\
\nonumber &=& \left(W^{1/2}B \left(x_{opt} - \tilde{x}_{opt}\right)\right)^T \left(W^{1/2}B \left(x_{opt} - \tilde{x}_{opt}\right)\right)\\
\label{eqn:energynorm} &=& \TNormS{W^{1/2}B \left(x_{opt} - \tilde{x}_{opt}\right)}.
\end{eqnarray}
We now use the formulas for $x_{opt}$ and $\tilde{x}_{opt}$, namely $x_{opt} = L^{\dagger}b$ (from eqn. (\ref{eqn:LaplacianOptVec})) and $\tilde{x}_{opt} = \tilde{L}^{\dagger}b$. Let $\Phi \in \mathbb{R}^{m \times n}$ denote the matrix $W^{1/2}B$ and let the SVD of $\Phi$ be
\begin{equation}\label{eqn:PhiSVD}
\Phi = U_{\Phi}\Sigma_{\Phi}V_{\Phi}^T.
\end{equation}
Here $U_\Phi \in \mathbb{R}^{m \times \rho}$, $\Sigma_\Phi \in \mathbb{R}^{\rho \times \rho}$, and $V_\Phi \in \mathbb{R}^{n \times \rho}$, with $\rho \leq n$ being the rank of $\Phi$.
Then, $$L = \Phi^T \Phi = V_{\Phi}\Sigma_{\Phi}^2 V_{\Phi}^T,$$ and thus
\begin{equation}\label{eqn:xoptSVD}
x_{opt} = L^{\dagger}b = V_{\Phi}\Sigma_{\Phi}^{-2} V_{\Phi}^T b.
\end{equation}
Similarly, $$\tilde{L} = \Phi^T SS^T \Phi = \left(S^{T} \Phi\right)^T \left(S^{T} \Phi\right)$$ and
\begin{equation}\label{eqn:xtildeoptSVD}
\tilde{x}_{opt} =  \left(S^{T} \Phi\right)^{\dagger} \left(S^{T} \Phi\right)^{\dagger T} b = \left(S^{T} U_{\Phi}\Sigma_{\Phi}V_{\Phi}^T\right)^{\dagger} \left(S^{T} U_{\Phi}\Sigma_{\Phi}V_{\Phi}^T\right)^{\dagger T}b .
\end{equation}
Combining eqns. (\ref{eqn:energynorm}), (\ref{eqn:PhiSVD}), (\ref{eqn:xoptSVD}), and (\ref{eqn:xtildeoptSVD}), we get
\begin{eqnarray}
\nonumber \left\|x_{opt} - \tilde{x}_{opt} \right\|_L  &=& \TNormS{U_{\Phi}\Sigma_{\Phi}V_{\Phi}^T \left(V_{\Phi}\Sigma_{\Phi}^{-2} V_{\Phi}^T b -
\left(S^{T} U_{\Phi}\Sigma_{\Phi}V_{\Phi}^T\right)^{\dagger} \left(S^{T} U_{\Phi}\Sigma_{\Phi}V_{\Phi}^T\right)^{\dagger T}b \right)}\\
\label{eqn:energynorm2} &=& \TNormS{\Sigma_{\Phi}^{-1} V_{\Phi}^T b -
\Sigma_{\Phi}\left(S^{T} U_{\Phi}\Sigma_{\Phi}\right)^{\dagger} \left(S^{T} U_{\Phi}\Sigma_{\Phi}\right)^{\dagger T}V_{\Phi}^Tb}.
\end{eqnarray}
In the above we used the facts that $U_{\Phi}$ and $V_{\Phi}$ are orthogonal matrices, and $\left(XV^T\right)^{\dagger} = VX^{\dagger}$ for any orthogonal matrix $V$. 
We now employ Theorem~\ref{thm:theorem7correct} of the Appendix in order to 
argue that $S^{T}U_{\Phi}$ is a matrix whose singular values are all close 
to unity. 
(This theorem is a variant of a result of Rudelson and Vershynin~\cite{RV07} 
that was proven as Theorem 4 in the appendix of~\cite{DMMS07_FastL2_TRv3}.) 
More specifically, since $U_{\Phi}^T U_{\Phi} = I_{\rho}$, Theorem~\ref{thm:theorem7correct} argues that with our choice of $r$ at Step 4 of Algorithm~\ref{alg:alg_sample_fast}
\begin{equation*}
\Expect{\TNorm{U_{\Phi}^T SS^T U_{\Phi}-I_{\rho}}} \leq \frac{\sqrt{\epsilon}}{6}.
\end{equation*}
Markov's inequality now implies that with probability at least $2/3$
\begin{equation}\label{eqn:rv}
\TNorm{U_{\Phi}^T SS^T U_{\Phi}-I_{\rho}} \leq \frac{\sqrt{\epsilon}}{2}.
\end{equation}
Using standard perturbation theory~\cite{Stewart90}, we get that for all $i=1,\ldots,\rho$,
\begin{equation}\label{eqn:svbound}
\abs{\sigma_i\left(U_{\Phi}^T SS^T U_{\Phi}\right) - 1} = \abs{\sigma_i^2\left(S^{T} U_{\Phi}\right)-1} \leq \frac{\sqrt{\epsilon}}{2}
\end{equation}
holds with probability at least 2/3. (Here $\sigma_i(X)$ denotes the $i$-th singular value of $X$.) This implies that the $m \times \rho$ matrix $S^{T}U_{\Phi}$ has rank $\rho$ with probability at least 2/3. The remainder of the proof will be conditioned on this event holding. Using $\left(S^TU_{\Phi}\Sigma_{\Phi}\right)^{\dagger} = \Sigma_{\phi}^{-1}\left(S^TU_{\Phi}\right)^{\dagger}$ (which is only true if $S^TU_{\Phi}$ has full rank), eqn. (\ref{eqn:energynorm2}) becomes
\begin{equation}
\label{eqn:energynorm3} \left\|x_{opt} - \tilde{x}_{opt} \right\|_L  = \TNormS{\Sigma_{\Phi}^{-1} V_{\Phi}^T b -
\left(S^{T} U_{\Phi}\right)^{\dagger} \left(S^{T} U_{\Phi}\right)^{\dagger T}\Sigma_{\Phi}^{-1}V_{\Phi}^Tb}.
\end{equation}
We now focus on the matrix $\Omega = S^{T}U_{\Phi} \in \mathbb{R}^{m \times \rho}$. Let its SVD be
\begin{equation}
S^{T}U_{\Phi} = \Omega = U_{\Omega} \Sigma_{\Omega} V_{\Omega}^T.
\end{equation}
Since the rank of $S^{T}U_{\Phi}$ is $\rho$, it follows that $U_\Omega \in \mathbb{R}^{m \times \rho}$, $\Sigma_\Omega \in \mathbb{R}^{\rho \times \rho}$, and $V_\Omega \in \mathbb{R}^{\rho \times \rho}$. We now rewrite eqn. (\ref{eqn:energynorm3}) using the SVD of $\Omega$:
\begin{equation}
\label{eqn:energynorm4} \left\|x_{opt} - \tilde{x}_{opt} \right\|_L  = \TNormS{\Sigma_{\Phi}^{-1} V_{\Phi}^T b -
V_{\Omega}\Sigma_{\Omega}^{-2}V_{\Omega}^T\Sigma_{\Phi}^{-1}V_{\Phi}^Tb}.
\end{equation}
Let $\Sigma_{\Omega}^{-2} = I_{\rho} + E$, for some diagonal error matrix $E$. Using $V_{\Omega}V_{\Omega}^T=V_{\Omega}^T V_{\Omega} = I_{\rho}$, eqn. (\ref{eqn:energynorm4}) becomes
\begin{eqnarray}
\nonumber \left\|x_{opt} - \tilde{x}_{opt} \right\|_L  &=& \TNormS{\Sigma_{\Phi}^{-1} V_{\Phi}^T b -
V_{\Omega}\left(I+E\right)V_{\Omega}^T\Sigma_{\Phi}^{-1}V_{\Phi}^Tb}\\
\nonumber &=& \TNormS{V_{\Omega}EV_{\Omega}^T\Sigma_{\Phi}^{-1}V_{\Phi}^Tb}\\
\nonumber &=& \TNormS{EV_{\Omega}^T\Sigma_{\Phi}^{-1}V_{\Phi}^Tb}\\
\nonumber &\leq& \TNormS{EV_{\Omega}^T}\TNormS{\Sigma_{\Phi}^{-1}V_{\Phi}^Tb}\\
\label{eqn:energynorm5} &=& \TNormS{E}\TNormS{\Sigma_{\Phi}^{-1}V_{\Phi}^Tb}.
\end{eqnarray}
We now seek to bound the spectral norm of the diagonal matrix $E$. Notice that the diagonal entries of $E$ satisfy
$$\abs{E_{ii}} = \abs{\sigma_{i}^{-2}\left(\Omega\right)-1} = \abs{\sigma_{i}^{-2}\left(S^{T}U_{\Phi}\right)-1}.$$
Using the bounds of eqn. (\ref{eqn:svbound}) we get
\begin{eqnarray}
\nonumber \TNorm{E} &=& \max_{i =1\ldots \rho} \abs{\sigma_{i}^{-2}\left(S^{T}U_{\Phi}\right)-1}\\
\nonumber &=& \max_{i =1\ldots \rho} \abs{\frac{\sigma_{i}^{2}\left(S^{T}U_{\Phi}\right)-1}{\sigma_{i}^{2}\left(S^{T}U_{\Phi}\right)}}\\
\label{eqn:energynorm6} &\leq& \frac{\sqrt{\epsilon}/2}{1-\left(\sqrt{\epsilon}/2\right)}\leq \sqrt{\epsilon}.
\end{eqnarray}
The last inequality follows since $\epsilon \leq 1$. Combining eqns. (\ref{eqn:energynorm5}) and (\ref{eqn:energynorm6}), we get
\begin{eqnarray}
\label{eqn:energynorm7} \left\|x_{opt} - \tilde{x}_{opt} \right\|_L &\leq& \epsilon\TNormS{\Sigma_{\Phi}^{-1}V_{\Phi}^Tb}.
\end{eqnarray}
To conclude the proof, notice that using $\Phi = W^{1/2}B$ and eqns. (\ref{eqn:PhiSVD}) and (\ref{eqn:xoptSVD}), we get
\begin{eqnarray}
\nonumber \left\|x_{opt}\right\|_L  &=& x_{opt}^T  L x_{opt}\\
\nonumber &=& \left(W^{1/2} B x_{opt}\right)^T  \left(W^{1/2} B x_{opt}\right) \\
\nonumber &=& \TNormS{\Phi x_{opt}}\\
\nonumber &=& \TNormS{U_{\Phi}\Sigma_{\Phi} V_{\Phi}^T V_{\Phi}\Sigma_{\Phi}^{-2} V_{\Phi}^T b}\\
\label{eqn:energynorm8}&=& \TNormS{\Sigma_{\Phi}^{-1} V_{\Phi}^T b}.
\end{eqnarray}
Combining eqns. (\ref{eqn:energynorm7}) and (\ref{eqn:energynorm8}) concludes the proof of the theorem.
\end{Proof}

\subsection{Running time}\label{sxn:runtime} 

We now discuss the running time of Algorithm~\ref{alg:alg_sample_fast}. 
Steps 1 and 2 are trivial and run in $O(m)$ time. 
Step 3 necessitates the computation of a probability distribution over the 
rows of $BW^{1/2}$. 
Theoretically, this step runs (for $\beta=1$) in $O(m \log^{c_1} n)$ time, 
for some small constant $c_1$, as described in~\cite{BSS09a_STOC}. 
(However, in order to achieve this running time it is necessary to perform 
$O(\log n)$ calls to the Spielman-Teng solver, which essentially renders 
this computation impractical. 
Below, we will discuss in more detail several issues related to computing
these probabilities in other ways.)
Steps 5, 6, and 7 run in $O(m)$ time, since $B$ is a matrix with two 
non-zero elements per row, $W$ is a diagonal matrix, and the sampling 
matrix $S$ simply reduces the number of rows in $BW^{1/2}$ from $m$ to $r$. 
Finally, at the last step, we invoke a direct solver for the sparse 
least-squares problem of eqn.~(\ref{eqn:LaplacianLS_Sample}), which takes 
$O\left(\frac{n^2}{\epsilon} \log\frac{n}{\epsilon}\right)$ time. 
Thus, from a theoretical perspective, using the fact that $m \leq n^2$, the 
running time Algorithm~\ref{alg:alg_sample_fast} is 
$O\left(\frac{n^2}{\epsilon} \left(\log\frac{n}{\epsilon}\right)\left(\log^{c_1} n\right)\right)$.


\section{Connecting graph resistances and statistical leverage scores}
\label{sxn:graph_resistances}

In this section, we will show that the effective resistances of the edges of a graph $G$ with $n$ vertices and $m$ positively weighted undirected edges are proportional to the statistical leverage scores of the rows of the matrix $W^{1/2}B$ (recall our definitions in Section~\ref{sxn:LaplacianDef}). Although this connection is straightforward from technical perspective, it is of considerable interest due to the insights it provides.

\subsection{Review of effective resistance and statistical leverage}\label{sxn:definitions}

We start with the following definition of the \emph{effective resistance} of an edge of a graph:
\begin{definition}
Given $G=(V,E)$, a connected, weighted, undirected graph with $n$ nodes, $m$ edges, and corresponding edge weights $w_e \ge 0$, for all $e\in E$, let
\begin{equation}
\label{eqn:LaplacianMatrix}
L = B^T W B
\end{equation}
denote the $n \times n$ Laplacian matrix of $G$ (see Section~\ref{sxn:LaplacianDef} for notation). The effective resistances $R_e$ across all edges $e \in E$ are
given by the diagonal entries of the matrix
\begin{equation}
\label{eqn:EffectiveResistance}
R= B L^{\dagger} B^T,
\end{equation}
where $L^{\dagger}$ denotes the Moore-Penrose generalized inverse of $L$.
\end{definition}
%
Clearly, from standard matrix algebra, the effective resistances of all the
edges of $G$ can be computed in $O(n^3)$ time.
Moreover, if we let $G$ denote an electrical network, in which each edge $e \in E$
corresponds to a resistor of resistance $1/w_e$, then the effective
resistance $R_e$ between two vertices can be defined as the potential
difference induced between the two vertices when a unit of current is
injected at one vertex and extracted at the other vertex.
Finally, effctive resistances have a wide range of applications, including not only theoretical applications such as analyzing diffusion processes and random walks on graphs, but also very practical applications such as analyzing clustering and community structure in large informatics networks.

A seemingly-unrelated notion is that of the \emph{statistical leverage
scores} of the rows of a matrix:
\begin{definition}
Given an $m \times n$ matrix $A$, with $m > n$, the \emph{statistical
leverage scores} of the rows of $A$ are the $m$ diagonal elements of the
projection matrix onto the span of the columns of~$A$.
That is, if the matrix $U_A$ denotes \emph{any} orthogonal basis for
the column space of $A$, then the diagonal elements of the projection
matrix $P_A$ onto the span of those columns are given by
$$
 (P_A)_{ii} = (U_{A}U_{A}^T)_{ii} = \VTTNormS{(U_{A})_{(i)}} ,
$$
where $(U_A)_{(i)}$ denotes the $i$-th row of the matrix $U_A$.
\end{definition}
%
Clearly, all the statistical leverage scores can be computed in $O(mn^2)$ time. Note that these scores could be defined for any $m \times n$ matrix $A$ with $m \le n$. In that case, however, if $A$ is not rank-deficient, then all the scores are trivially equal to unity.
Importantly, the statistical leverage scores have a natural interpretation in terms of ``importance'' or ``influence'' or ``leverage'' of the corresponding constraint/row of~$A$ in the overconstrained least squares optimization problem $ \min_x \TNorm{Ax-b}$. As such, they have been of interest historically in diagnostic regression
analysis~\cite{ChatterjeeHadi88}.

More generally, given a rank parameter $k$, one can define the \emph{statistical leverage scores relative to the best rank-$k$ approximation to $A$} to be the $m$ diagonal elements of the projection
matrix onto the span of the best rank-$k$ approximation to $A$. These generalized scores have been used recently as importance sampling probabilities to obtain relative-error approximation algorithms for regression~\cite{DMM06,DMMS07_FastL2_TRv3}, and they were essential for for the extension of these ideas to relative-error low-rank matrix approximation~\cite{DMM08_CURtheory_JRNL,CUR_PNAS} problems.
Prior work~\cite{CR07,AMT09_DRAFT} has also used term \emph{incoherent} to 
refer to the situation when no leverage score is particularly large.

\subsection{A simple lemma}\label{sxn:connection}

We now describe a connection between graph resistances and statistical leverage scores. Although this connection is not so surprising from a technical perspective---indeed, it is obvious once it is pointed out---it is useful for the insights it provides.

\begin{lemma}
Let the matrix $\Phi = W^{1/2}B \in \mathbb{R}^{m \times n}$ denote the edge-incidence matrix of a graph $G$ rescaled by $W^{1/2}$. The \emph{statistical leverage scores} associated with $\Phi$ are
(up to scaling) equal to the \emph{effective resistances} of all edges of a weighted graph $G$. That is, if $\ell_i$ is the leverage score associated with the $i$-th row of $\Phi$, then $\ell_i/w_i$ is the effective resistance of the $i$-th edge.
\end{lemma}
\begin{Proof}
Consider the matrix
\begin{equation*}
\label{eqn:Pmatrix}
P = W^{1/2}B(B^TWB)^+B^TW^{1/2} \in \mathbb{R}^{m \times m}  ,
\end{equation*}
and notice that $P=W^{1/2}RW^{1/2}$ is simply a rescaled version of the $m \times m$ matrix $R = B L^+ B^T$, whose diagonal entries are exactly equal to the effective resistances of all the edges of $G$. Since $\Phi = W^{1/2}B$,
it follows that
\begin{equation*}
 P = \Phi(\Phi^T\Phi)^+\Phi^T.
\end{equation*}
Let $U_{\Phi}$ denote an orthogonal matrix spanning the column space of $\Phi$. Then $ P = U_{\Phi}U_{\Phi}^T $, from which it follows that the
diagonal elements of $P$ are equal to
\begin{equation*}
\label{eqn:Pmatrix_diagonals}
P_{ii} = (U_{\Phi}U_{\Phi}^T)_{ii}
       = \VTTNormS{(U_{\Phi})_{(i)}}.
\end{equation*}
This concludes the proof of the lemma.
\end{Proof}

\subsection{Usefulness of statistical leverage in randomized matrix algorithms}
\label{sxn:graph_resistances:usefulness}

The connection between statistical leverage and effective resistance is of 
interest in attempts to make nearly-linear-time linear equation solvers more 
practical.
The reason is that statistical leverage has proven to be the key structural
quantity to understand in order to bridge the ``theory-practice gap'' between
theoretical work on randomized algorithms for large matrices, and
applications (both numerical-implementation and data-analysis applications)
of this ``randomized matrix algorithm''
paradigm~\cite{Paschou07b,CUR_PNAS,DMMS07_FastL2_TRv3,RT08,AMT09_DRAFT}.
In this section, we review some of the ``lessons learned,'' in the hope that
they provide insights on how to to bridge the theory-practice gap for
solving linear equations defined by a Laplacian constraint matrices.

Recall that much work, including, \emph{e.g.}, our previous
work~\cite{dkm_matrix1,dkm_matrix2,dkm_matrix3}, followed that of Frieze,
Kannan, and Vempala~\cite{FKV04}, in which columns and/or rows from a
matrix $A$ are randomly sampled according to a probability distribution
that depends on the Euclidean norms of those columns/rows.
In this case, worst-case additive-error guarantees of the form
\begin{equation}
\FNorm{A-P_{C,k}A} \le \FNorm{A-A_k} + \epsilon \FNorm{A}
\label{eqn:cx_additive}
\end{equation}
can be obtained, with high probability.%
\footnote{Here $P_{C,k} A$ denotes the projection of $A$ on a rank-$k$ subspace spanned by the columns of $C$.}
Although these algorithms were motivated by resource-constrained
computational environments, they have several drawbacks with respect to
numerical applications and data analysis applications more generally.
First, worst-case additive-error bounds are quite coarse.
Second, these algorithms were not immediately-relevant to common problems, as
they are typically formulated, in scientific computing and numerical linear
algebra.
Third, the insights provided by the sampling probabilities into the data
are limited---the probabilities are often uniform due to data preprocessing,
or they may correspond, e.g., simply to the degree of a node if the
data matrix is derived from a graph.

Importantly, each
of these three problems was solved by the introduction of importance
sampling probabilities that depend on the statistical leverage scores.%
\footnote{Although these probabilities were introduced
in~\cite{DMM06,DMM08_CURtheory_JRNL} and were used in solving two
very traditional numerical linear algebra problems
in~\cite{DMMS07_FastL2_TRv3,BMD09_CSSP_SODA}, the connection with leverage
scores wasn't made explicit until~\cite{CUR_PNAS}.}
\begin{itemize}
\item
First, by using importance sampling probabilities that depend on the
leverage scores, it was shown~\cite{DMM08_CURtheory_JRNL,CUR_PNAS} that
one could randomly sample a ``small'' number of columns to obtain
worst-case relative-error guarantees of the form
\begin{equation}
\FNorm{A-P_{C,k}A} \le (1+\epsilon) \FNorm{A-A_k}   ,
\label{eqn:cx_relative}
\end{equation}
with high probability.
\item
Second, algorithms that were comparable to or better than
previously-existing algorithms were provided for the following two very
traditional scientific computing problems:
\begin{itemize}
\item
\textbf{Overconstrained Least Squares.}
Let $A$ be an $m \times n$ matrix $A$, with $m \gg n$, and consider solving
$ x_{opt} = \arg \min_x \TNorm{Ax - b} $.
In previous work~\cite{DMM06,DMM08_CURtheory_JRNL,DMMS07_FastL2_TRv3},
we proposed a simple, sampling-based, algorithm for solving this problem:
first, compute the statistical leverage scores of the rows of $A$; then,
use these scores to construct an importance sampling probability
distribution to sample a ``small'' number of rows of $A$ and the
corresponding elements of $b$; and finally, solve the induced, much smaller
but still overconstrained, regression problem using only those (suitably
rescaled) rows of $A$ and the corresponding elements of $b$.
Strong relative error guarantees for this overconstrained%
\footnote{Note that it is easy to show that similar results hold for the
very underconstrained problem.
Let $A$ be an $m \times n$ matrix, with $m \ll n$, and consider the problem
of finding the minimum-length solution to
$ x_{opt} = \argmin_{x}||Ax-b||_2 = A^+b$.
Sampling variables or columns from $A$ can be represented by postmultiplying
$A$ by a $n \times c$ (with $c>m$) column-sampling matrix $S$ to
construct the (still underconstrained) least-squares problem:
$ \tilde{x}_{opt} = \argmin_{x}||ASS^Tx-b||_2 = A^T(AS)^{T+}(AS)^{+}b $.
The second equality follows by inserting $P_{A^T}=A^TA^{T+}$ to obtain
$ASS^TA^TA^{T+}x-b$ inside the $||\cdot||_2$ and recalling that
$A^+=A^TA^{T+}A^+$ for the Moore-Penrose pseudoinverse.
If one randomly samples $c=O((n/\epsilon^2) \log(n /\epsilon))$ columns
according to ``column-leverage-score'' probabilities, i.e., the
diagonal elements of the projection matrix onto the row space, then it can
be proven that $||x_{opt}-\tilde{x}_{opt}||_2 \le \epsilon||x_{opt}||_2$
holds, with high probability.}
regression problem were proven with this
approach~\cite{DMM06,DMM08_CURtheory_JRNL}.
\item
\textbf{Column Subset Selection Problem.}
Let $A$ be an $m \times n$ matrix, and let $k$ be a positive
integer.
Then, pick $k$ columns of $A$ forming an $m \times k$ matrix $C$ such that
the residual $\XNorm{A - P_C A}$, where $\xi = 2\ \mbox{or}\ F$ denotes the
spectral norm or Frobenius norm, is minimized over all possible
${n \choose k}$ choices for the matrix $C$.
Previously \cite{BMD09_CSSP_SODA,BMD08_CSSP_KDD}, we developed a two-phase
algorithm that uses the nonuniformity structure defined by the statsitical
leverage scores in an essential way to provide theoretical and empirical
results for both the spectral and Frobenius norm that were competitive
or better than previously existing results.
\end{itemize}
\item
Third, the insights into the matrix provided by statistical leverage scores
(in both numerical and data applications) can be quite refined.
The insights are used in very different ways, depending on whether one is
interested in high-quality numerical implementations or large-scale data
analysis applications.
\begin{itemize}
\item
\textbf{Numerical Implementation Applications.}
Here, one wants to provide fast high-quality numerical implementations, and
one is typically interested in the error parameter to be bery small,
\emph{e.g.}, $\epsilon \approx 10^{-16}$.
For example, with respect to the overconstrained least-squares regression
problem, performing an exact computation of the statistical
leverage scores of the rows of $A$ is no faster than exactly solving the
original regression problem.
Sarl\'{o}s~\cite{Sarlos06,DMMS07_FastL2_TRv3} addressed this problem by
preprocessing the matrix $A$ and the vector $b$ with the randomized Hadamard
transform of Ailon and Chazelle~\cite{AC06}.
This preprocessing step made the statistical leverage scores almost
uniform---effectively ``washing out'' any nonuniformities defined by the
leverage scores, thereby densifying the matrix if it was sparse---thus
leading to the first randomized, relative-error algorithm for least-squares
problems that runs asymptotically faster than $\Theta(mn^2)$ time.
High-quality implementations of such algorithms have
appeared~\cite{RT08,AMT09_DRAFT}, and they highlight the significant
practical applicability of this approach.
\item
\textbf{Data Analysis Applications.}
Here, one may want $\epsilon \approx 0.1$, and one is typically interested
in obtaining insight with respect to some downstream data analysis goal.
In such cases, SVD-based methods are often chosen for computational
convenience, rather than because the statistical assumptions underlying
their use are satisfied by the data---a fact which means that the leverage
scores are often extremely nonuniform in a way that correlates strongly
with what practitioners know about the
data~\cite{Paschou07b,CUR_PNAS,BMD08_CSSP_KDD,BMD09_kmeans_NIPS} problems.
Thus, far from ``washing out'' this nonuniformity structure, one is
interested in identifying and exploiting it.
Intuitively, conditioned on
being reliable, more ``outlier-like'' data points may be the most important
and informative.
\end{itemize}
\end{itemize}

This brings us to the question of how to compute these statistical leverage
scores, or equivalently the effective resistances, which is an issue that 
gets to the heart of the theory-practice gap.
Depending on the application and the resource constraints, there are several 
alternatives:
\begin{itemize}
\item
Compute the scores by calling the Spielman-Teng nearly-linear time solver. 
This algorithm runs in $O\left(\nnz{A} \log^{c_1} n\right)$ time, where 
$\nnz{A}$ represents the number of non-zero elements of the matrix $A$, or 
equivalently the number of edges in the graph $G$, and $c_1$ is a small 
constant.
This method works for computing the leverage scores of Laplacian matrices; 
and in this case it is, theoretically, the best method.
\item
Compute the scores by computing an ``exact'' basis for the column space of
the $m \times n$ matrix $\Phi = W^{1/2}B$.
This takes $O(mn^2)$ time and works for general matrices.
For Laplacian matrices it is clearly expensive, given that the 
weighted edge-incidence matrix is very sparse.
\item
Compute an approximation to the scores based on iterative sampling and 
volume sampling ideas that have been used in relative-error low-rank matrix 
approximations~\cite{DRVW06_JRNL,DR10_TR}.
This might be of interest if a pass-efficient model is an appropriate model 
for data access.
\item
Compute an approximation to the scores based on numerical methods to, 
e.g., compute an estimator for the diagonal of a matrix~\cite{BKS07}.
These numerical methods are particularly appropriate for large matrices when 
matrix-vector products are easy to evaluate; they have proven useful in 
uncertainty quantification~\cite{BCF10}; and they draw on the observation 
that the leverage scores, being proportional to the diagonal elements of a 
projection matrix, have a natural interpretation in scientific computing in 
terms of density matrices and Green's functions~\cite{SCS10}.
\end{itemize}
These alternate approaches are of particular interest since data points with 
high leverage scores often have natural interpretations in terms of 
processes generating the data matrices~\cite{CUR_PNAS}.
Moreover, an examination of the details of these methods illustrates that 
problems are parameterized within theoretical computer science in very 
different ways than they are parameterized in scientific computing.
Finally, an important issue to keep in mind is that in most applications, 
one does not need a uniformly good approximation to all the leverage scores, 
but instead one needs a good approximation only to the ``high leverage'' 
data points.


\section{Conclusion}
\label{sxn:discussion}

Several open problems suggest themselves.
On the theoretical side: 
Can one draw on the original ideas of Spielman and Teng in order to
develop an algorithm with the simplicity of ours and with the running time
approximation of theirs?
Similarly, can we get the $O(n \log n)$ factor, which currently is due to 
the result of Rudelson and Vershynin~\cite{RV07}, down to $O(n)$, even for 
some classes of graphs, thereby obtaining a more immediately practical 
version of the result of Batson, Spielman, and Srivastava~\cite{BSS09a_STOC}?
On the more applied side:
How rapidly can we approximate (even with a one-sided approximation)
the statistical leverage scores, either for general $m \times n$ matrices
$A$ and arbitrary rank parameter $k$, or under some realistic generative
model?
Similarly, can one use the connection between statistical leverage and 
effective resistance to design improved heuristics, given knowledge about
the processes generating the data?

We conclude by noting that the last two questions are of particular interest.
Although much of the recent work on using Laplacian preconditioners has
focused on nearly-linear-time solvers for computing ``exact'' solutions,
\emph{i.e.}, with the error parameter $\epsilon$ set to machine precision,
there are many other applications of these ideas.
For example, in machine learning, Ravikumar and Lafferty used preconditioner
approximations for doing approximate inference in probabilistic graphical
models~\cite{RL06}.
This connection should not be surprising, as much of the work on the
``randomized algorithms for matrices'' paradigm has been motivated by
large-scale data applications.
In many of these data analysis applications, however, not only is setting
$\epsilon=10^{-16}$ not of interest, doing so would actually lead to
``worse'' answers than setting it, say, as $\epsilon=0.1$.
If other recent applications of the randomized algorithms paradigm are any
guide~\cite{RT08,AMT09_DRAFT,CUR_PNAS,BW09_PNAS}, then the issues that will
arise when thinking of $\epsilon$ as extremely small and trying to couple
newer randomized algorithmic methods with traditional numerical
methods~\cite{RT08,AMT09_DRAFT} will be very different than the issues that
arise in applications where the data are much less well-structured and
much-coarser $\epsilon$'s are of interest~\cite{CUR_PNAS,BW09_PNAS}.


\section*{Appendix}

Let $A \in \mathbb{R}^{m \times n}$ be any matrix. 
Consider the following algorithm, which is essentially the algorithm in 
page 876 of~\cite{DMM08_CURtheory_JRNL}.
This algorithm constructs a matrix $C\in \mathbb{R}^{m \times c}$ consisting 
of $c$ sampled and rescaled columns of $A$. 

\begin{algorithm}[ht]
\begin{framed}

\SetLine

\AlgData{
$A \in \mathbb{R}^{m \times n}$,
$p_i \geq 0, i\in[n]$ s.t. $\sum_{i \in [n]}p_i=1$,
positive integer $c \leq n$.}

\AlgResult{
$C \in \mathbb{R}^{m \times c}$
}

Initialize $S \in \mathbb{R}^{m \times c}$ to be an all-zero matrix.

\For{$t=1,\ldots,c$}{
   Pick $i_t \in [n]$, where $\textbf{Prob}\left(i_t = i\right) = p_i$\;
   $S_{i_t t} = 1/\sqrt{cp_{i_t}}$\;
}

Return $C = AS$\;

\end{framed}
\caption{
The \textsc{Exactly($c$)} algorithm.
}
\label{alg:SDconstruct_exact}
\end{algorithm}

\noindent
Next, we state a theorem that provides a bound for the approximation error 
$\TNorm{AA^T-CC^T}$.
We used this in the proof of our main theorem in 
Section~\ref{sxn:LinearEqnSolving} in order to argue that the singular 
values of the ``sampled orthogonal'' matrix $S^{T}U_{\Phi}$ are all close 
to unity. 
In this form, the theorem was proven as Theorem 4 in the Appendix 
of~\cite{DMMS07_FastL2_TRv3}, but it is a variant of the well-known result 
of Rudelson and Vershynin~\cite{RV07}. 

\begin{theorem}\label{thm:theorem7correct}
Let $A \in \mathbb{R}^{m \times n}$ with $\TNorm{A} \leq 1$. Construct $C$ using the \textsc{Exactly($c$)} algorithm and let the sampling probabilities $p_i$ satisfy
\begin{equation}\label{eqn:defPj}
p_i \geq \beta \frac{\TNormS{A^{(i)}}}{\FNormS{A}}
\end{equation}
for all $i \in [n]$ for some constant $\beta \in (0,1]$. Let $\epsilon \in (0,1)$ be an accuracy parameter, assume $c_0^2 \FNormS{A}\geq 4\beta \epsilon^2$, and let
\begin{equation*}\label{eqn:Cbound}
c = 2 \left(\frac{c_0^2 \FNormS{A}}{\beta \epsilon^2}\right)\log \left(\frac{c_0^2 \FNormS{A}}{\beta \epsilon^2}\right).
\end{equation*}
(Here $c_0$ is the unknown constant of Theorem 3.1, p. 8 of~\cite{RV07}.) Then,
$$\Expect{\TNorm{AA^T-CC^T}}\leq \epsilon.$$
\end{theorem}
Finally, it is worth noting that the condition 
$c_0^2 \FNormS{A}\geq 4\beta \epsilon^2$ is trivially satisfied for any 
matrix $A$ such that $\FNormS{A} \geq 4$ assuming $c_0 \geq 1$.

\end{document}